\documentclass[12pt]{article}
\usepackage{amssymb}

\newcommand{\RR}{{\mathbb R}}            % the reals R
\begin{document}
\begin{center}
\Large
% TITLE GOES HERE
Uncountable sets and an infinite real number game
\end{center}

\begin{flushright}
Matthew H.~Baker
\footnote{Supported by NSF Grant DMS-0300784.} \\
School of Mathematics \\
Georgia Institute of Technology \\
Atlanta, GA 30332-0160 \\
\verb+mbaker@math.gatech.edu+
\end{flushright}

% \begin{abstract}
% We give a short proof of the well-known fact that the unit interval $[0,1]$ in $\RR$ is uncountable 
% by means of a simple infinite game.   We also show using this game that a (non-empty) perfect subset of $[0,1]$
% must be uncountable.  
% \end{abstract}

\paragraph*{The game.}

Alice and Bob decide to play the following infinite game on the real number line.
A subset $S$ of the unit interval $[0,1]$ is fixed, and then Alice and Bob alternate playing real numbers.
Alice moves first, choosing any real number $a_1$ strictly between $0$ and $1$.  Bob then chooses any
real number $b_1$ strictly between $a_1$ and $1$.  On each subsequent turn, the players must choose a point
strictly between the previous two choices.  Equivalently, if we let $a_0 = 0$ and $b_0 = 1$, then 
in round $n\geq 1$, Alice chooses a real number $a_n$ with $a_{n-1} < a_n < b_{n-1}$, and then 
Bob chooses a real number $b_n$ with $a_n < b_n < b_{n-1}$.
Since a monotonically increasing sequence of real numbers which is bounded
above has a limit (see \cite[Theorem 3.14]{Rudin}), 
$\alpha = \lim_{n \to \infty} a_n$ is a well-defined real number between 0 and 1.
Alice wins the game if $\alpha \in S$, and Bob wins if $\alpha \not\in S$.

\paragraph*{Countable and uncountable sets.}

An set $X$ is called {\em countable} if it is possible to list the elements of $X$ in a
(possibly repeating) infinite sequence $x_1,x_2,x_3,\ldots$
% \footnote{According to our definition, finite sets are countable.  Some authors require }
Equivalently, $X$ is countable if there is a function from the 
set $\{ 1,2,3, \ldots \}$ of natural numbers to $X$ which is {\em onto}.
For example, every finite set is countable, and the set of natural numbers is countable.
A set which is not countable is called {\em uncountable}.  Cantor proved using his famous {\em diagonalization argument} 
that the real interval $[0,1]$ is uncountable.  We will give a different proof of this fact based on Alice and Bob's game.

\medskip

\textsc{Proposition 1.}
{\it If $S$ is countable, then Bob has a winning strategy.}

{\it Proof.}
Since $S$ is countable, one can enumerate the elements of $S$ as $s_1,s_2,s_3,\ldots$
Consider the following strategy for Bob.  On move $n\geq 1$, he chooses $b_n = s_n$ if this is a legal move,
and otherwise he randomly chooses any allowable number for $b_n$.  Since $\alpha < b_n$ for all $n$, it follows that
$\alpha \neq b_n$ for any $n \geq 1$, and thus $\alpha \not\in S$.  This means that Bob always wins with this strategy!

\medskip

If $S = [0,1]$, then clearly Alice wins no matter what either player does.  Therefore we deduce:

\textsc{Corollary 1.}
{\it The interval $[0,1] \subset \RR$ is uncountable.}

\medskip

This argument is in many ways much simpler than Cantor's original proof!

\paragraph*{Perfect sets.}

We now prove a generalization of the fact that $[0,1]$ is uncountable.  This will also follow
from an analysis of our game, but the analysis is somewhat more complicated.
Given a subset $X$ of $[0,1]$, we make the following definitions:

\begin{itemize}
\item A {\em limit point} of $X$ is a point $x \in [0,1]$ such that for every $\epsilon > 0$,
the open interval $(x - \epsilon , x + \epsilon)$ contains an element of $X$ other than $x$.
\item $X$ is {\em closed} if every limit point of $X$ belongs to $X$.
\item $X$ is {\em perfect} if it is non-empty
\footnote{Some authors consider the empty set to be perfect.}, closed, and if 
every element of $X$ is a limit point of $X$.  
\end{itemize}

For example, the famous middle-third {\em Cantor set} is perfect (see \cite[\S2.44]{Rudin}).
If $L(X)$ denotes the set of limit points of $X$, then a nonempty set $X$ is 
closed $\Leftrightarrow L(X) \subseteq X$, and is perfect
$\Leftrightarrow L(X) = X$.
% From the definition, one sees that a non-empty set $X$ is perfect precisely when 
% the set of limit points of $X$ is $X$ itself.
It is a well-known fact that every perfect set is
uncountable (see \cite[Theorem 2.43]{Rudin}).
Using our infinite game, we will give a different proof of this fact.
We recall the following basic property of the interval $[0,1]$:

\begin{itemize}
\item[($\star$)] Every non-empty subset $X \subseteq [0,1]$ has an {\em infimum} (or {\em greatest lower bound}),
meaning that
there exists a real number $\gamma \in [0,1]$ such that
$\gamma \leq x$ for every $x \in X$, and if $\gamma' \in [0,1]$ is any real number with
$\gamma' \leq x$ for every $x \in X$, then $\gamma' \leq \gamma$.  
\end{itemize}

The infimum $\gamma$ of $X$ is denoted by 
${\rm inf}(X)$.
% \item[(ii)] For any two real numbers $a,b \in [0,1]$ with $a < b$, there exists a real number
% $c$ with $a < c < b$.
% (This easily implies that for any real numbers $a,b \in [0,1]$ with $a < b$, there exist {\em infinitely many} real numbers
% $c$ with $a < c < b$.)

\medskip

Let's say that a point $x \in [0,1]$ is {\em approachable from the right}, denoted $x \in X^+$,
if for every $\epsilon > 0$, the open interval $(x , x + \epsilon)$ contains an element of $X$.
We can define {\em approachable from the left} (written $x \in X^-$) similarly using the interval $(x - \epsilon, x)$.
It is easy to see that $L(X) = X^+ \cup X^-$, so that a non-empty set $X$ is perfect $\Leftrightarrow X = X^+ \cup X^-$.

\medskip

The following two lemmas tell us about approachability in perfect sets.

\medskip

\textsc{Lemma 1.}
{\it If $S$ is perfect, then $\inf(S) \in S^+$.}

{\it Proof.}
The definition of the infimum in ($\star$) implies that ${\rm inf}(S)$ cannot be approachable from the left,
so, being a limit point of $S$, it must be approachable from the right.

\medskip

\textsc{Lemma 2.}
{\it If $S$ is perfect and
$a \in S^+$, then for any $\epsilon > 0$, the
open interval $(a,a+\epsilon)$ also contains an element of $S^+$.}

{\it Proof.}
Since $a \in S^+$, we can choose three points $x,y,z \in S$ with $a < x < y < z < a + \epsilon$.
Since $(x,z) \cap S$ contains $y$, the real number $\gamma = \inf((x,z) \cap S)$ 
satisfies $x \leq \gamma \leq y$.
If $\gamma = x$, then by $(\star)$ we have $\gamma \in S^+$.  If $\gamma > x$, then $(\star)$ implies that $\gamma \in L(X)$ and $(x,\gamma) \cap S = \emptyset$, so that
$\gamma \not\in S^-$ and therefore $\gamma \in S^+$.

\medskip

From these lemmas, we deduce:

\textsc{Proposition 2.}
{\it If $S$ is perfect, then Alice has a winning strategy.}

{\it Proof.}
Alice's only constraint on her $n$th move is that $a_{n-1} < a_n < b_{n-1}$.
By induction, it follows from Lemmas 1 and 2
that Alice can always choose $a_n$ to be an element of $S^+ \subseteq S$.
Since $S$ is closed, $\alpha = \lim a_n \in S$, so Alice wins!

\medskip

From Propositions 1 and 2, we deduce:

\textsc{Corollary 2.}
{\it Every perfect set is uncountable.}

\paragraph*{Further analysis of the game.}

We know from Proposition 1 that Bob has a winning strategy if $S$ is countable, and
it follows from Proposition 2 that Alice has a winning strategy if $S$ contains a perfect set.  (Alice just
chooses all of her numbers from the perfect subset.)  What can one say in general?  
A well-known result from set theory \cite[\S6.2, Exercise 5]{Ciesielski} says that every uncountable {\em Borel set}
\footnote{A Borel set is, roughly speaking, any subset of $[0,1]$ that can be constructed by taking countably many unions, intersections, and complements of open intervals;
see \cite[\S11.11]{Rudin} for a formal definition.}
contains a perfect subset.
Thus we have completely analyzed the game when $S$ is a Borel set: Alice wins if $S$ is uncountable, and Bob wins
if $S$ is countable.
However, there do exist non-Borel 
% (and therefore difficult to describe)
uncountable subsets of $[0,1]$ which do not contain a perfect subset \cite[Theorem 6.3.7]{Ciesielski}.  
So we leave the reader with the following problem:

\medskip

{\bf Problem:} Do there exist uncountable subsets of $[0,1]$ for which:
(a) Bob has a winning strategy; 
(b) Alice does not have a winning strategy; or
(c) neither Alice nor Bob has a winning strategy?

% \begin{itemize}
% \item[(a)] Bob has a winning strategy;
% \item[(b)] Alice does not have a winning strategy; or
% \item[(c)] neither Alice nor Bob has a winning strategy?
% \end{itemize}

\paragraph*{Related games.}
% \paragraph*{Origins of the game.}

Our infinite game is a slight variant of the one proposed by Jerrold Grossman and Barry Turett in \cite{GT}
(see also \cite{Newcomb}). Propositions 1 and 2 above
were motivated by parts (a) and (c), respectively, of their problem.
The author originally posed Propositions 1 and 2
as challenge problems for the students in his 
Math 25 class at Harvard University in Fall 2000.
%, but didn't write up the solution until $5\frac{1}{2}$ years later!

\medskip

A related game (the ``Choquet game'') can be used to prove the Baire category theorem (see
\S8.C of \cite{Kechris} and \cite{HL}).  In Choquet's game, played in a given metric space $X$, Pierre moves first by choosing a non-empty
open set $U_1 \subseteq X$.  Then Paul moves by choosing a non-empty open set $V_1 \subseteq U_1$.
Pierre then chooses a non-empty open set $U_2 \subseteq V_1$, etc., yielding two decreasing sequences
$U_n$ and $V_n$ of non-empty open sets with $U_n \supseteq V_n \supseteq U_{n+1}$ for all $n$, and $\cap U_n = 
\cap V_n$.  Pierre wins if $\cap U_n = \emptyset$, and Paul wins if $\cap U_n \neq \emptyset$.  One can show that
if $X$ is complete, then Paul has a winning strategy, 
and if $X$ contains a non-empty open set $O$ which is a countable union of closed sets
having empty interior, then Pierre has a winning strategy.  As a consequence, one obtains the {\em Baire category theorem:} 
If $X$ is a complete metric space, then no open subset of $X$ can be a countable union of closed sets having empty interior.

\medskip

Another related game is the Banach-Mazur game (see \S6 of \cite{Oxtoby} and \S8.H of \cite{Kechris}).
A subset $S$ of the unit interval $[0,1]$ is fixed, and then Anna and Bartek alternate play.
First Anna chooses a closed interval $I_1 \subseteq [0,1]$, and then Bartek chooses a closed interval $I_2 \subseteq I_1$.
Next, Anna chooses a closed interval $I_3 \subseteq I_2$, and so on.  Together the players' moves determine a nested sequence
$I_n$ of closed intervals.  Anna wins if $\cap I_n$ has at least one point in common with $S$, otherwise Bartek wins.
% Banach (unpublished) proved Mazur's conjecture 
It can be shown that Bartek has a winning strategy if and only if $S$ is meagre
(see Theorem 6.1 of \cite{Oxtoby}).  
(A subset of $X$ is called {\em nowhere dense} if the interior of its closure is empty, and
is called {\em meagre}, or of the {\em first category}, if it is a countable union of nowhere dense sets.)
It can also be shown, using the axiom of choice,
that there exist sets $S$ for which the Banach-Mazur game is undetermined (neither player has a winning strategy).

\medskip

For a more thorough discussion of these and many other {\em topological games}, we refer the reader to the
survey article \cite{Telgarsky}, which contains an extensive bibliography.  Many of the games discussed in
\cite{Telgarsky} are not yet completely understood.

\medskip

Games like the ones we have been discussing play a prominent role in the modern field of {\em descriptive set theory}, most notably in
connection with the {\em axiom of determinacy} (AD).
(See Chapter 6 of \cite{Kanamori} for a more detailed discussion.)
Let $X$ be a given subset of the space $\omega^\omega$ of infinite sequences of natural numbers, and consider the following game
between Alice and Bob.  Alice begins by playing a natural number, then Bob plays another (possibly the same) natural number, then
Alice again plays a natural number, and so on.  The resulting sequence of moves determines an element $x \in \omega^\omega$.  Alice wins
if $x \in X$, and Bob wins otherwise.  
% A theorem of Martin shows that this game is determined (i.e., one of the players has a winning strategy)
% when $X$ is a Borel set.  
The axiom of determinacy states that this game is determined 
(i.e., one of the players has a winning strategy) for {\em every} choice of $X$.

\medskip

A simple construction shows that the axiom of determinacy is inconsistent with the axiom of choice.
On the other hand, with Zermelo-Fraenkel set theory plus the axiom of determinacy (ZF+AD), 
one can prove many non-trivial theorems about the real numbers, including: (i) every subset of $\RR$ is
Lebesgue measurable; and (ii) every uncountable subset of $\RR$ contains a perfect subset.
Although ZF+AD is not considered a ``realistic'' alternative to ZFC (Zermelo-Fraenkel + axiom of choice), it has stimulated
a lot of mathematical research, and certain variants of AD are taken rather seriously.  
For example, the axiom of {\em projective determinacy} is intimately connected with the continuum hypothesis and
the existence of large cardinals (see \cite{Woodin} for details).

% \newpage


\begin{thebibliography}{10}

\bibitem{Ciesielski} K.~Ciesielski, {\em Set Theory for the Working Mathematician}, London Mathematical Society Student Texts {\bf 39}, 
Cambridge University Press, 1997.

\bibitem{GT} J.~W.~Grossman and B.~Turett, Problem \#1542, {\em Mathematics Magazine} {\bf 71}, February 1998.

\bibitem{HL} F.~Hirsch and G.~Lacombe, {\em Elements of Functional Analysis}, Graduate Texts in Mathematics {\bf 192},
Springer-Verlag, 1999.

\bibitem{Kanamori} A. Kanamori, {\em The Higher Infinite} (2nd ed.), Springer-Verlag, 2003.

\bibitem{Kechris} A. Kechris, {\em Classical Descriptive Set Theory}, Springer-Verlag, 1995.

\bibitem{Newcomb} W.~A.~Newcomb, Solution to Problem \#1542, {\em Mathematics Magazine} {\bf 72}, February 1999.

\bibitem{Oxtoby} J. Oxtoby, {\em Measure and Category} (2nd ed.), Springer-Verlag, 1980.

\bibitem{Rudin} W. Rudin,
{\em Principles of Mathematical Analysis} (3rd ed.), McGraw-Hill, 1976.

\bibitem{Telgarsky} R. Telg{\'a}rsky, {\em Topological games: On the 50th anniversary of the Banach-Mazur game}, Rocky Mountain J. Math. {\bf 17} (1987), 
227--276.

\bibitem{Woodin} H. Woodin, {\em The Continuum Hypothesis, Part I}, Notices of the AMS {\bf 48}, no. 6 (2001), 567--576.

\end{thebibliography}
\end{document}